\documentclass[final]{siamltex}

\usepackage{amsfonts}
\usepackage{amssymb,amsmath}

\newcommand{\ud}{\mathrm{d}}
\newcommand{\e}{\mathrm{e}}

\begin{document}

\title{An \emph{a posteriori} condition on the numerical approximations of the Navier-Stokes equations for the existence of a strong solution\footnote{\scriptsize{Submitted to SIAM Journal on Numerical Analysis}}}
\author{ Masoumeh Dashti \and James C. Robinson \thanks{Mathematics Institute, University of Wawick, Coventry CV4 7AL. UK ({\tt m.dashti@warwick.ac.uk}), ({\tt j.c.robinson@warwick.ac.uk}).}}

\maketitle

\begin{abstract}
In their 2006 paper, Chernyshenko et al prove that a sufficiently smooth strong solution of the 3d Navier-Stokes equations is robust with respect to small enough changes in initial conditions and forcing function. 
They also show that if a regular enough strong solution exists then Galerkin approximations converge to it. They then use these results to conclude that the existence of a sufficiently regular strong solution can be verified using sufficiently refined numerical computations. 
In this paper we study the solutions with minimal required regularity to be strong, which are less regular than those considered in Chernyshenko et al (2006). We prove a similar robustness result and show the validity of the results relating convergent numerical computations and the existence of the strong solutions.

\end{abstract}

\begin{keywords} 
Navier-Stokes equations, Galerkin method
\end{keywords}

\section{Introduction}
For the three-dimensional Navier-Stokes equations the existence of strong solutions over
 an arbitrary time interval is not known. Therefore the validity of the results of the numerical solutions of 3d Navier-Stokes equations is not obvious. This problem is addressed in Chernyshenko et al (2006) where a rigorous relationship between numerical and sufficiently regular exact solutions is given. They show for sufficiently smooth initial conditions and forcing functions (data) that although \emph{a priori} there is no guarantee of the validity of the numerical solutions, there is an \emph{a posteriori} condition that if satisfied by the numerical results guarantees the existence of a strong solution. In this paper we show the validity of the results proved in Chernyshenko et al (2006) for the less regular strong solutions which are not covered there.\\
\indent We will study the Navier-Stokes equations in their functional form. For a bounded domain $\Omega$ we let $\mathcal{H}$ be the space of divergence-free smooth vector-valued functions on $\Omega$ with compact support and zero average and define 
\begin{eqnarray*}
H &=& \textrm{closure of } \mathcal{H}\: \mathrm{in} \:[L^2(\Omega)]^3,\\
V &=& \textrm{closure of } \mathcal{H}\: \mathrm{in} \:[H^1(\Omega)]^3.
\end{eqnarray*}
We use the same notation $H$ and $V$ for the similar spaces of periodic functions over the periodic domain $Q$. Then the Navier-Stokes equations in their functional form are written as (Constantin and Foias 1988, Robinson 2001)
\begin{equation}\label{NSe_reg}
\frac{\ud u}{\ud t} +\nu Au +B(u,u) = f,\quad\mbox{with}\quad u(0)=u_0
\end{equation}
where $Au=-\Pi\Delta$, $B(u,u)=\Pi (u.\nabla)u$ with $\Pi$ the orthogonal projection from $L^2$ into $H$.
We will consider this equation with the following cases for the data\\
\\
\indent (a) `minimal regularity' when $u_0\in{V}$ and $f\in L^2(0,T;H)\cap L^1(0,T;V)$\\
\indent (b) `second order regularity' when $u_0\in{V^2}$ and $f\in L^2(0,T;V)\cap L^1(0,T;V^2)$,\\
\\
where $V^m=H^m\cap V$. For the  periodic case we know that $D(A^{m/2})=H^m\cap V$ for all $m$ and therefore we define the norm on $V^m$ as $\|u\|_m=|A^{m/2}u|$. We denote by $(\cdot,\cdot)$ and $|\cdot|$ the inner product and norm on $H$.\\
\indent We will show that the strong solution of (\ref{NSe_reg}) with the data introduced in $(a)$ or $(b)$, if it exists, remains strong if the changes in the initial condition and forcing function are small enough. The exact conditions required for these changes to be `small' are given in theorems \ref{robust1} and \ref{robust2}. For example in the case of minimally regular strong solutions we require
\begin{eqnarray*}
&&|D(u_0-v_0)| +\int_0^T|Df(s)-Dg(s)|\,\ud s \nonumber\\
&& \quad < \frac{1}{k}\left(\frac{\nu^3}{27T}\right)^{1/4}
                      \exp\left(-\frac{k^2}{2}\int_0^T\frac{27k^2}{2}\frac{1}{\nu^3}|Du(s)|^4+\frac{1}{\nu}|Du(s)|\,|Au(s)|\,\ud s\right).
\end{eqnarray*}
We then, in theorems \ref{postt} and \ref{postt2}, use these robustness results to find an \emph{a posteriori} condition that if satisfied by sufficiently refined numerical approximations, implies the existence of a strong solution. We also show that if a strong solution exists the Galerkin approximations converge to it and then use this to prove that the existence of a strong solution can be verified by the Galerkin approximations. In the last section we consider a channel flow as a physical example that can be described by the Navier-Stokes equations with the conditions introduced above. For this example we will show how the results of this paper can be applied to the Galerkin approximations to verify the existence of a strong solution.\\ 
\indent The results we prove here for a strong solution with lowest regularity hold in a general bounded domain as well as in the absence of boundaries unlike the results for more regular solutions which are proved only for the equations in a periodic domain or the whole space. 

%%%%%%%%%%%%%%%%%%%%%%%%%%%%%%%%%%%%%%%%%%%%%
%                                     ODE lemma
%%%%%%%%%%%%%%%%%%%%%%%%%%%%%%%%%%%%%%%%%%%%%
\section{General ODE lemma}
We first prove an ODE lemma which will be used in dealing with the differential inequalities that appear in the proofs. We consider the differential inequality
$$
\frac{\ud y}{\ud t}\le \delta(t)+\alpha y^n \qquad \mbox{with}\quad y(0)=y_0>0
$$
and find the conditions on $y_0$, $\delta(t)$ and $\alpha$ that ensure that $y(t)$ exists on a finite time interval $[0,T]$. This lemma is a generalization of the result obtained for $n=2$ in Chernyshenko et al (2006).
\begin{lemma}\label{lemma} %%%%%%%%%%% ODE Lemma
Let $T>0$, $\alpha>0$ and $n>1$ be constants and let $\delta(t)$
be a non-negative continuous function on $[0,T]$. Let $y$ satisfy
the differential inequality
\begin{equation} \label{ineq}
\frac{\ud y}{\ud t}\le \delta(t)+\alpha y^n \qquad \mbox{with}\quad y(0)=y_0>0
\end{equation}
and define
$$
\eta=y_0+\int_0^T\delta(s)\,\ud s.
$$
\begin{itemize}
\item[(i)] If
\begin{equation} \label{general_condition}
\eta<\frac{1}{[(n-1)\alpha T]^{1/(n-1)}}
\end{equation}
then $y(t)$ remains bounded on $[0,T]$
\item[(ii)]
$y(t)\rightarrow0$ uniformly on $[0,T]$ as $\eta\rightarrow0$.
\end{itemize}
\end{lemma}

\begin{proof}

We first consider the following differential inequality
\begin{equation}\label{simple}
\dot{z} \le \alpha z^n \quad\mbox{with}\quad z(0)=\eta
\end{equation}
and show that
\begin{equation}\label{sup_y(T)}
\sup_{S_1}\;{y(T)} \le \sup_{S_2}\;{z(T)}.
\end{equation}
where $S_1$ and $S_2$ are the sets of all possible solutions of inequalities (\ref{ineq}) and (\ref{simple}) respectively.\\
\indent Since $\dot{y}$, $\dot{z}$, $y(0)$ and $z(0)$ are non-negative, the suprema of $y(T)$ and $z(T)$ happen when 
\begin{eqnarray*}
\dot{y}&&=\delta(t)+ \alpha y^n,\\
\dot{z}&&=\alpha z^n
\end{eqnarray*}
for all $t\in [0,T]$. In this case, for the difference $w=y-z$ we have
\begin{equation}
\dot{w}=\delta(t)+ \alpha w \sum_{k=0}^{n-1} {n-1 \choose k} y^{n-1-k}\;z^k \quad \mbox{with}\quad w(0)=-\int_0^t\delta(s)\;\ud s.\nonumber
\end{equation}
Since $y(t)$ and $z(t)$ are greater than zero and assuming they remain finite over $[0,T]$, there exists some $M>0$ such that 
$$
\dot{w} \le \alpha M w + \delta(t)
$$
and therefore by Gronwall's inequality
\begin{eqnarray*}
w(t) &&\le w(0)\e^{\alpha Mt}+\int_0^t\delta(s)\;\ud s\\
     &&=(1-\e^{\alpha Mt})\int_0^t\delta(s)\;\ud s.
\end{eqnarray*}
This implies $w(t)\le 0$ for any $t\in[0,T]$ and (\ref{sup_y(T)}) follows.\\

The inequality (\ref{simple}) can be written as
$$
\frac{\dot{z}}{z^n}\le \alpha
$$
Integrating both sides from 0 to $t$ yields
\begin{equation}\label{z(T)}
z(T) \le \frac{\eta}{(-\alpha T (n-1)(\eta)^{n-1}+1)^{1/{n-1}}}.
\end{equation}
Since $y(t)\le y(T) \le z(T)$, $y(t)$ remains finite on $[0,T]$ provided that
$$
\alpha T (n-1)(\eta)^{n-1} < 1
$$
which yields (\ref{general_condition}).\\
\indent Furthermore, it is clear from (\ref{z(T)}) that $z(T)\to 0$ as $\eta\to 0$, from which it follows that $y(t)\to 0$ uniformly on $[0,T]$.
\end{proof}

%%%%%%%%%%%%%%%%%%%%%%%%%%%%%%%%%%%% NSE
%
\section{The inequalities for the nonlinear term}
\indent For the nonlinear operator, we will use the following inequalities in this paper
\begin{eqnarray}
|(B(u,v),Aw)|   &&\le k|Du|\,|Dv|^{1/2}|Av|^{1/2}|Aw|\label{triform1}\\
|(B(w,u),A^2 w)|&&\le c\|u\|_3\|w\|_2^2\label{triform21}\\ 
|(B(u,w),A^2 w)|&&\le c'\|u\|_3\|w\|_2^2\label{triform22},
\end{eqnarray}
The first inequality holds for both Dirichlet and periodic boundary conditions (Constantin and Foias 1988) while the other two are valid only in the absence of boundaries with $c$ and $c'$ independent of the size of the domain (Kato 1972). 
The constant $k$ for the periodic and no-slip domain and also the whole of $\mathbb{R}^3$, is independent of the domain $\Omega$. In a general domain with nonzero boundary conditions however, it depends on the regularity properties of the boundary of $\Omega$. To see this we use the H\"older inequality to write
$$
|B(u,v,Aw)| \le \sum_{i,j=1}^3 \|u_i\|_{L^6}\|D_iv_j\|_{L^3}|Aw_j|.
$$
By the Sobolev inequality
\begin{equation} \label{sobolev1}
\|u\|_{L^6}\le c_s|Du|
\end{equation}
we have
$$
\|Dv\|_{L^3}\le \|Dv\|_{L^6}^{1/2}|Dv|^{1/2}\le c_s^{1/2}|Av|^{1/2}|Dv|^{1/2}.
$$
Therefore we can write
\begin{eqnarray*}
|B(u,v,Aw)| &\le& \sum_{i,j=1}^3 \|u_i\|_{L^6}\|D_iv_j\|_{L^3}|Aw_j|\\
                  &\le& 9\|u\|_{L^6}\|Dv\|_{L^3}|Aw|\\
                  &\le& 9c_s^{3/2}|Du|\,|Dv|^{1/2}|Av|^{1/2}|Aw|
\end{eqnarray*}
and so $k=9c_s^{3/2}$. For no-slip boundary conditions and the whole of $\mathbb{R}^3$, $c_s$, the constant of the Sobolev inequality (\ref{sobolev1}), does not depend on $\Omega$ (Ziemer 1989).
For a bounded domain with non-zero boundary conditions, $c_s$ depends on the regularity properties of $\partial\Omega$ but is independent of the size of $\Omega$ (Adams 1975). Therefore for the cubic domain of the periodic boundary conditions also, it does not depend on the size of the domain. In fact, the cubic domain has the strong Lipschitz property and for such a domain Adams (1975, Lemma 5.10) has shown that $c_s=4\sqrt{2}$. 
Similarly, $c$ and $c'$ depend on the constant of the Sobolev inequality $\|u\|_{L^{6/(3-2k)}}\le c_{s,k}\|u\|_{H^k}$ which is again independent of the size of the domain (Adams 1975 and Ziemer 1989).\\
\indent 
In obtaining the robustness conditions in the next section, we keep track of the constants $k$, $c$ and $c'$ when they appear, so that the value of the constant coefficients in the robustness conditions could be computed if desired.\\
\indent The inequalities (\ref{triform1})--(\ref{triform22}) are not as elegant as the inequalities available for the more regular solutions considered in Chernyshenko et al (2006) and this is why different proofs are needed in the less regular cases we are studying here.\\
\indent We note here that assuming a strong solution $u\in L^\infty(0,T;V)\cap L^2(0,T;V^2)$ exists, from (\ref{triform1}) we have $B(u,u)\in L^2(0,T;H)$ and therefore (\ref{NSe_reg}) implies that $\ud u/\ud t \in L^2(0,T;H)$. Having $u_0\in{V^2}$ and $f\in L^2(0,T;V)$ in case (b), from the regularity result for periodic domains proved in Constantin and Foias (1988) we know that the strong solution $u$ is in fact more regular and an element of $L^\infty(0,T;V^2)\cap L^2(0,T;V^3)$.

%%%%%%%%%%%%%%%%%%%%%%%%%%%%%%%%%%%%%%%%%%%%%
%                               (Robustness)
%%%%%%%%%%%%%%%%%%%%%%%%%%%%%%%%%%%%%%%%%%%%%
\section{Robustness of strong solutions}

Here we show that as long as a strong solution exists for some specific initial data and forcing function, the equations with sufficiently close data also have a strong solution. 
A similar result about the robustness of strong solutions with respect to changes in the forcing function is proved by Fursikov in his 1980 paper. For the three-dimensional Navier-Stokes equations with initial condition $u_0\in V^m$ where $m\ge 1/2$, he shows that the set of forcing functions for which a strong solution exists is open in $L^q(0,T;V^{m-1})$ with $q\ge 2$.
\\
\indent For both minimally regular solutions (a) and the more regular case (b) the condition we obtain here for the data depends explicitly on the viscosity coefficient unlike the more regular cases considered in Chernyshenko et al (2006). The robustness result we prove for the minimally regular strong solutions holds in a general bounded domain. For the second order regular solutions however, we need to restrict the domain to be periodic or the whole of $\mathbb{R}^3$ since the inequalities (\ref{triform21}) and (\ref{triform22}) only hold in a periodic domain (Constantin and Foias 1988) or the whole of $\mathbb{R}^3$ (Kato 1972).
\subsection*{(a) Strong solutions with minimal regularity}
To prove the robustness with respect to the data we write the governing equation for the difference between the solution of the equations with nearby data and the strong solution and find a bound on the norm of this difference. 
\begin{theorem}\label{robust1} %%%%%%%% Robustness of strong solutions for m=1
Let $f\in L^1(0,T,V)\cap L^2(0,T;H)$, $u_0\in V$ and $u\in L^{\infty}(0,T;V)\cap L^2(0,T;V^2)$ be a strong solution of
\[
\frac{\ud u}{\ud t}+Au+B(u,u)=f \quad\mbox{with}\quad u(0)=u_0.
\]
If
\begin{eqnarray}\label{m1cond}
&&|D(u_0-v_0)| +\int_0^T|Df(s)-Dg(s)|\,\ud s \nonumber\\
&& \quad < \frac{1}{k}\left(\frac{\nu^3}{27T}\right)^{1/4}
                      \exp\left(-\frac{k^2}{2}\int_0^T\frac{27k^2}{2}\frac{1}{\nu^3}|Du(s)|^4+\frac{1}{\nu}|Du(s)|\,|Au(s)|\,\ud s\right)
\end{eqnarray}
then the solution of
\[
\frac{\ud v}{\ud t}+\nu Av+B(v,v)=g\quad\mbox{with}\quad v(0)=v_0
\]
is also a strong solution on $[0,T]$ with the same regularity as $u$.
\end{theorem}

\begin{proof}
By the local existence of strong solutions (Constantin and Foias 1988, Temam 1995) we know that there exists $T^*>0$ such that $v\in L^{\infty}(0,T';V)\cap L^2(0,T';V^2)$ for any $T'<T^*$. We consider $T^*$ the maximal time of existence of the strong solution $v$ meaning that $\lim\sup_{t\to T^*} |Dv|=\infty$. In the following argument we assume $T^*\le T$ and deduce a contradiction.\\
\indent We consider $w=u-v$ which satisfies
$$
\frac{\ud w}{\ud t}+\nu Aw+B(u,w)+B(w,u)-B(w,w)=f-g.
$$
Over the time interval $t\in [0,T')$ for any $T'<T^*$ we have $dv/dt\in L^2(0,T';H)$ and since $T^*\le T$, $du/dt\in L^2(0,T';H)$. Therefore taking the inner product of the above equation with $Aw$ and using (\ref{triform1}) we can write
\begin{eqnarray*}
\frac{1}{2}\frac{\ud}{\ud t}|Dw|^2 + \nu|Aw|^2 
               &\le&  k|Du||Dw|^{1/2}|Aw|^{3/2}+ k|Dw||Du|^{1/2}|Au|^{1/2}|Aw|\\
               &    &  \qquad+k|Dw|^{3/2}|Aw|^{3/2} +|D(f-g)|\;|Dw|.
\end{eqnarray*}
We then use Young's inequality to remove $|Aw|$ (which causes the appearence of $\nu$ in the coefficients) and obtain
$$
\frac{1}{2}\frac{\ud}{\ud t}|Dw|^2 \le \left( \frac{27k^4}{4\nu^3} |Du|^4 +\frac{k^2}{2\nu}|Du||Au| \right)|Dw|^2+ \frac{27k^4}{4\nu^3}|Dw|^6 +|D(f-g)|\;|Dw|.
$$
Dividing by $|Dw|$ (assuming $|Dw|$ is nonzero. If not, replacing $|Dw|$ by $|Dw|+e_0$ with $e_0>0$ will lead us to the same result after tending $e_0$ to zero at the end) we get
$$
\frac{\ud}{\ud t}|Dw| \le \left( \frac{27k^4}{4\nu^3} |Du|^4 +\frac{k^2}{2\nu}|Du||Au| \right)|Dw|+ \frac{27k^4}{4\nu^3}|Dw|^5 +|D(f-g)|.
$$
Now letting
$$
\beta(t)= \frac{k^2}{2} \int_0^t  \left( \frac{27k^2}{2}\frac{1}{\nu^3}|Du|^4 +\frac{1}{\nu}|Du|\,|Au| \right)\; \ud s
$$
and setting $y(t)=|Dw(t)|\e^{-\beta(t)}$ the above inequality can be written as
$$
\frac{\ud y}{\ud t} \le \alpha y^5 + \delta(t),
$$
where $\alpha=\frac{27k^4}{4\nu^3}\e^{4\beta(T)}$ and $\delta(t)=|Df(t)-Dg(t)|$.\\
\indent By lemma \ref{lemma}, if the condition (\ref{m1cond}) is satisfied, $|Dw(t)|$ is uniformly bounded on
 $[0,T^*)$. This implies that $|Dv(T^*)|<\infty$ which contradicts the maximality of $T^*$. It follows
 that $v(t)$ is a strong solution on $[0,T]$.
\end{proof} %------------------------------------------------------------
\subsection*{(b) Strong solutions with second order regularity}

We follow the same approach as in the minimally regular case to show the robustness for a strong solution with second order regularity. Here also the inequalities we need to use for the nonlinear operator to find a bound on the $V^2$-norm of the difference between nearby solutions depend on the $V^3$-norm of this difference and therefore the robustness condition on the data depends explicitly on the viscosity coefficient. However, its dependence is weaker in this case.
\begin{theorem}\label{robust2} %%%%%% Robustness of strong solutions for m=2
Let $f\in L^1(0,T,V^2)\cap L^2(0,T;V)$, $u_0\in V^2$ and $u\in L^{\infty}(0,T;V^2)\cap L^2(0,T;V^3)$ be a
 strong solution of
$$
\frac{\ud u}{\ud t}+Au+B(u,u)=f \quad\mbox{with}\quad u(0)=u_0.
$$
If
\begin{eqnarray}\label{m2cond}
|A(u_0-v_0)| +\int_0^T|A(f-g)|\,\ud t\;   <  \; \frac{1}{c}\sqrt{\frac{2\nu}{T}}\exp\left(-\int_0^T (c+c')\|u\|_3\, \ud t\right),
\end{eqnarray}
then the solution of
$$
\frac{\ud v}{\ud t}+Av+B(v,v)=g \quad\mbox{with}\quad v(0)=v_0
$$
is also a strong solution on $[0,T]$ with the same regularity as $u$.
\end{theorem}  %%%%%%%%%%%%%%%%%%%%%%%%%%%%%%%%%%%%%%

\begin{proof}
As before, the difference $w=u-v$ satisfies
\begin{equation}
\frac{\ud w}{\ud t}+\nu Aw+B(u,w)+B(w,u)+B(w,w)=f-g.\nonumber
\end{equation}
Taking the inner product of the above equation with $A^2 w$ and using  (\ref{triform21}) and (\ref{triform22}) we obtain
\begin{equation}
\frac{1}{2}\frac{\ud}{\ud t}\|w\|_2+ \nu\|w\|_3^2 \le (c+c')\|u\|_3^2\|w\|_2^2 + c\|w\|_3\|w\|_2^2 + \|f-g\|_2\|w\|_2.\nonumber
\end{equation}
We apply Young's inequality to the second term on the right hand side and get
$$
\frac{1}{2}\frac{\ud}{\ud t}\|w\|_2^2 \le (c+c')\|u\|_3\|w\|_2^2 + \frac{c^2}{4\nu}\|w\|_2^4 +\|f-g\|_2\|w\|_2.
$$
Dividing the above inequality by $\|w\|_2$ yields
$$
\frac{\ud}{\ud t}\|w\|_2 \le (c+c')\|u\|_3\, \|w\|_2 + \frac{c^2}{4\nu}\|w\|_2^3 +\|f-g\|_2.
$$
\indent Now we consider
\[
y(t)=\|w\|_2 e^{-\beta(t)}
\]
with $\beta(t)=\int_0^t (c+c')\|u\|_3\, \ud t$. The equation with this new variable becomes
\[
\frac{\ud y}{\ud t}\le \delta(t)+\alpha y^3,
\]
where $\alpha=\frac{c^2}{4\nu}e^{2\beta(T)}$ and $\delta=\|f-g\|_2$.
Using lemma \ref{lemma} $\|w\|_2$ is bounded if the condition (\ref{m2cond}) is satisfied.
\end{proof} %------------------------------------------------------------------
%

%%%%%%%%%%%%%%%%%%%%%%%%%%%%%%%%%%%%%%%%%%%%%
%                               Convergence of GA
%%%%%%%%%%%%%%%%%%%%%%%%%%%%%%%%%%%%%%%%%%%%%

\section{Convergence of Galerkin approximations}

Here we show that if a strong solution with minimal or second order regularity exists, Galerkin approximations converge to it. Similar results about convergence of various numerical method assuming the existence of a strong solution, are given for finite element methods by Heywood and Rannacher (1982), for a Fourier collocation method by E (1993) and for a nonlinear Galerkin method by Devulder, Marion and Titi (1993). Here as in Chernyshenko et al (2006), we make no assumption on the regularity of the Galerkin approximations.\\
\indent For the minimally regular case the result we prove here holds for the solution of the equations over a general bounded domain. For the second order regular solutions we need to use inequality (\ref{triform22}) which is valid only in a periodic domain or the whole of $\mathbb{R}^3$.\\ 
\indent In the following theorems we let $P_n$ be the orthogonal projection in $H$ onto the space spanned by the first $n$ eigenfunctions of the Stokes operator $A$, $\{w_j\}_{j=1}^n$, ordered so that their corresponding eigenvalues satisfy 
$0 < \lambda_1 \le \lambda_2 \le \dots$ . 
Since the eigenfunctions of $A$ are smooth (Constantin and Foias 1988), for any $u\in V^m$, $m\ge 0$ (with $V^0=H$) we have $u_n=P_n u=\sum_{j=1}^n (u,w_j)w_j\in V^m$.\\
\indent We note that what we obtain here about the convergence of Galerkin approximations would be useful even if the existence of regular solutions was known.

\subsection*{(a) Minimal regularity} 

The following theorem, like the robustness theorem in the minimal regularity case, holds in sufficiently smooth bounded domains as well as in the absence of boundaries.
\begin{theorem}\label{galerkin1} %%%%%%%%%% m=1 Galerkin Convergence
Let $u_0\in V$, $f\in L^2(0,T;H)$ and $u\in L^{\infty}(0,T;V) \cap L^2(0,T;V^2)$ be a strong solution of the Navier-Stokes equations
$$
\frac{\ud u}{\ud t} +\nu Au +B(u,u)=f(t)\quad\mbox{with}\quad u(0)=u_0.
$$
Then $u_n$, the solution of Galerkin system
\begin{equation}\label{galerkin_sys}
\frac{\ud u_n}{\ud t} +\nu Au_n +P_n B(u_n,u_n)=P_n f(t)\quad\mbox{with}\quad u_n(0)=P_n u_0,
\end{equation}
converges strongly to $u$ in both $L^{\infty}(0,T;V)$ and $L^2(0,T;V^2)$ as $n\to \infty$.
\end{theorem}
\begin{proof}
We consider $w_n=u-u_n$ which satisfies
$$
\frac{\ud w_n}{\ud t} +\nu Aw_n +P_n B(u,w_n) +P_n B(w_n,u) -P_n B(w_n,w_n) =Q_n f(t) -Q_n B(u,u).
$$
Letting $q_n=u-P_n u$ and taking the inner product of the above equation with $Aw_n$ while noting that
 $(P_n B(u,v),Aw_n)=b(u,v,P_nAw_n)$ we obtain
\begin{eqnarray*}
\frac{1}{2}\frac{\ud}{\ud t}|Dw_n|^2 +\nu|Aw_n|^2 
                          \le &&\,k\left(|Du|\,|Dw_n|^{1/2}|Aw_n|^{3/2}\right. \\
                               &&+|Dw_n|\,|Du|^{1/2}\,|Au|^{1/2}|Aw_n|
                                                               +|Dw_n|^{3/2}|Aw_n|^{3/2}\\
                               &&\left.+|Du|^{3/2}|Au|^{1/2}|Aq_n|  \right)+|f|\,|Aq_n|.
\end{eqnarray*}
Applying Young's inequality and redefining the constant $k$ gives
\begin{eqnarray*}
\frac{\ud}{\ud t}|Dw_n|^2 \le &&\,k\left\{\left( \frac{1}{\nu^3}|Du|^4    
                                                  +\frac{1}{\nu}|Du|\,|Au|\right)|Dw_n|^2 
                                                  +\frac{1}{\nu^3}|Dw_n|^6\right.\\
                                            &&\left.+|Du|^{3/2}|Au|^{1/2}|Aq_n|\right\}+|f|\,|Aq_n|.
\end{eqnarray*}
Now letting 
$$
\beta(t)=k\int_0^t\left( \frac{1}{\nu^3}|Du|^4 +\frac{1}{\nu}|Du|\,|Au|\right)\,\ud s
$$ 
and 
$$
y_n(t)=\mathrm{e}^{-\beta(t)}|Dw_n|^2
$$
the above inequality can be written as
$$
\frac{\ud y_n}{\ud t} \le \alpha y_n^3 +\delta_n(t),
$$
with $\alpha=\frac{k}{\nu^3}\mathrm{e}^{2\beta(T)}$ and $\delta_n(t)=|f|\,|Aq_n| +k|Du|^{3/2}|Au|^{1/2}|Aq_n|$.
\\
\indent Using the H\"older inequality we have
\begin{eqnarray*}
\int_0^T \delta_n(s)\,\ud s \le \left[\left(\int_0^T |f|^2\ud s\right)^{1/2}
                               +k\left(\int_0^T |Du|^{3}|Au|\ud s\right)^{1/2}\right] 
                                          \left(\int_0^T |Aq_n|^2\ud s\right)^{1/2}
\end{eqnarray*}
Since $u\in L^2(0,T;V^2)$, $u(s)\in V^2$ for almost every $s\in [0,T]$ and therefore
 $q_n(s)=Q_n u(s)\to 0$ in $V^2$ as $n\to\infty$ and $|Aq_n(s)|^2\le |Au(s)|^2$ for a.e. $s\in [0,T]$. Hence the
 Lebesgue dominated convergence theorem for each $s$ implies that $\int_0^T |Aq_n|^2\ud s$ tends to zero
 as $n\to \infty$ and therefore 
$$
\int_0^T \delta_n(s)\,\ud s\to 0
$$
as $n\to \infty$. Since $y_n(0)\to 0$ as $n\to\infty$, lemma \ref{lemma} shows that
$y_n$, and hence $|Dw_n|^2$, converges to zero uniformly on
$[0,T]$ as $n\to\infty$.

\end{proof} %------------------------------------------------------------

\subsection*{(b) Second order regularity}

Here also, in a similar way to the previous case, we try to find a bound on the $V^2$-norm of the difference between a second order regular strong solution and Galerkin approximations. However here the more regular solution space makes the proof easier.
\begin{theorem}\label{galerkin2} %%%%%%%%%% m=2 Galerkin Convergence
Let $u_0\in V^2$, $f\in L^2(0,T;V)\cap L^1(0,T;V^2)$ and $u\in L^{\infty}(0,T;V^2) \cap L^2(0,T;V^3)$ be 
a strong solution of the Navier-Stokes equations
$$
\frac{\ud u}{\ud t} +\nu Au +B(u,u)=f(t)\quad\mbox{with}\quad u(0)=u_0.
$$
Then $u_n$, the solution of Galerkin system
\begin{equation}\label{galerkin_sys2}
\frac{\ud u_n}{\ud t} +\nu Au_n +P_n B(u_n,u_n)=P_n f(t)\quad\mbox{with}\quad u_n(0)=P_n u_0,
\end{equation}
converges strongly to $u$ in both $L^{\infty}(0,T;V^2)$ and $L^2(0,T;V^3)$ as $n\to \infty$.
\end{theorem}
\begin{proof}
Let $w_n=u-u_n$. Then $w_n$ satisfies
$$
\frac{\ud w_n}{\ud t} +\nu Aw_n +P_n B(u,w_n) +P_n B(w_n,u) +P_n B(w_n,w_n) = Q_n f -Q_n B(u,u)
$$
Taking the inner product of the above equation with $A^2 w_n$ and using (\ref{triform21}) and
 (\ref{triform22}) we obtain
\begin{eqnarray*}
\frac{1}{2}\frac{\ud}{\ud t}\|w_n\|_2^2 +\nu\|w_n\|_3^2  
                              &\le& (c+c')\|u\|_3\|w_n\|_2^2 +c\|w_n\|_3\|w_n\|_2^2\\
                              &   & +\|Q_nf-Q_n B(u,u)\|_2\|w_n\|_2
\end{eqnarray*}
Now we use Young's inequality to remove the dependence on $\|w_n\|_3$ and then divide by $\|w\|_2$ to obtain
$$
\frac{\ud}{\ud t}\|w_n\|_2 \le (c+c')\|u\|_3 \|w_n\|_2 +\frac{c^2}{4\nu}\|w_n\|_2^3 +\|Q_n f-Q_nB(u,u)\|_2,
$$
in which the coefficient of $\|w_n\|_2$ does not depend on $n$.\\
\indent Letting
$$
y_n=\|w_n\|_2 e^{-\beta(t)},
$$
with $\beta(t)=\exp\left( \int_0^t (c+c')\|u(s)\|_3\;\ud s\right)$ yields
$$
\dot{y}_n \le \delta(t) + \alpha y_n^3,
$$
where $\alpha=\frac{c^2}{4\nu}e^{2\beta(T)}$ and $\delta(t)=\|Q_n f-Q_nB(u,u)\|_2$.
So by lemma \ref{lemma} if 
$$
\eta = \|Q_n u_0\|_m + \int_0^T \|Q_n f-Q_nB(u,u)\|_2\;\ud s \to 0 \quad \mathrm{as}
 \quad n\to\infty,
$$
then $y_n(t)\to 0$ uniformly on $[0,T]$ as $n\to \infty$.\\
By (\ref{galerkin_sys2}) $\|Q_n u_0\|_m\to 0$ as $n\to\infty$. We know that 
\begin{equation}\label{triform23}
\|B(u,u)\|_2 \le c\|u\|_2\|u\|_3
\end{equation}
(Kato 1972, Constantin and Foias 1988) and therefore $f(s) -B(u(s),u(s)) \in V^2$ for almost every $s\in [0,T]$. So since $\{w_j\}_{j=1}^n$ form a basis for $V^2$ as well as $H$, $\|Q_n \left(f(s)-B(u(s),u(s))\right)\|_2\to 0$ and
$$
\|Q_n \left(f(s)-B(u(s),u(s))\right)\|_2 \le \|f(s)-B(u(s),u(s))\|_2
$$
for almost every $s\in [0,T]$. Therefore by the Lebesgue dominated convergence theorem it follows that 
$$
\int_0^T \|Q_n \left(f(s)-B(u(s),u(s))\right)\|_2 \to 0
$$ 
and the result follows.

\end{proof} %--------------------------------------------------------------------
%
%%%%%%%%%%%%%%%%%%%%%%%%%%%%%%%%%%%%%%%%%%%
%                               LINK
%%%%%%%%%%%%%%%%%%%%%%%%%%%%%%%%%%%%%%%%%%%
%
\section{Numerical verification of the existence of a strong solution}
Here we show that the existence of minimal and second order regular strong solutions can be verified via computations using sufficiently refined Galerkin approximations.
\begin{theorem}\label{postt}
$\mathbf{i)}$ Consider the Navier-Stokes equations 
\begin{equation}\label{post_nse}
\frac{\ud u}{\ud t}+\nu Au+B(u,u)=f\quad\mbox{with}\quad u(0)=u_0
\end{equation}
with $u_0\in V$, $f\in L^2(0,T;H)\cap L^1(0,T;V)$, that hold in a bounded domain $\Omega$ with sufficiently smooth boundary or in a periodic domain. \\
\indent Let $v\in L^{\infty}(0,T;V)\cap L^2(0,T;V^2)$ be a numerical approximation of $u$ satisfying
$$
\frac{\ud v}{\ud t}+\nu Av+B(v,v) \in L^1(0,T;V)\cap L^2(0,T;H)
$$
and
\begin{eqnarray}
&&|Dv(0)-Du_0| + \int_0^T \|\frac{dv(s)}{ds}+\nu Av(s)+B(v(s),v(s))
                                                                           -f(s)\|_1\,\ud s \nonumber\\
&&\quad < \frac{1}{k}\left(\frac{\nu^3}{27T}\right)^{1/4}
                      \exp\left(-\frac{k^2}{2}\int_0^T\frac{27k^2}{2}\frac{1}{\nu^3}|Dv(s)|^4
                      +\frac{1}{\nu}|Dv(s)||Av(s)|\,\ud s\right).\label{posteriori}
\end{eqnarray}
Then the Navier-Stokes equations (\ref{post_nse}) have a strong solution
 $u\in L^{\infty}(0,T;V)\cap L^2(0,T;V^2)$.\\ 
\\
\indent $\mathbf{ii)}$ Let $u$ be a strong solution of (\ref{post_nse}). Then there exists an $N$ such that the Galerkin approximations $u_n$ satisfy the inequality (\ref{posteriori}) for all $n>N$. Therefore in view of part $\mathbf{(i)}$, $u$ passes the  \emph{a posteriori} test as a solution approximated by $u_n$ i.e. the existence of a strong solution with minimal regularity can be verified by the Galerkin approximations. 
\end{theorem}
\begin{proof}

$\mathbf{i)}$ Considering
$$
g=\frac{\ud v}{\ud t}+\nu Av +B(v,v),
$$
$v$ is a strong solution of
$$
\frac{\ud \bar v}{\ud t}+\nu A\bar v +B(\bar v,\bar v)=g.
$$
So by theorem \ref{robust1}, if the inequality (\ref{posteriori}) holds, the solution of the equations with nearby data $(u_0,f)$ is a strong solution.\\
\\
\indent $\mathbf{ii)}$ The strong convergence of $u_n$ to $u$ in $L^{\infty}(0,T;V)$ and $L^2(0,T;V^2)$ is guaranteed by theorem \ref{galerkin1}. Therefore the right hand side of (\ref{posteriori}) is bounded below for every $n$. It follows from (\ref{galerkin_sys}) that $|Du_n(0)-Du_0|\to 0$ as $n\to 0$ and 
$$
\frac{du_n}{dt}+\nu Au_n +B(u_n,u_n)-f=Q_n[B(u_n,u_n)-f].
$$
Therefore it only remains to show $\|Q_n[B(u_n(s),u_n(s))-f]\|_1$ converges to zero as $n\to \infty$ for a.e. $s\in[0,T]$, it is then clear that (\ref{posteriori}) will be satisfied for all $n$ sufficiently large.\\
\indent For the nonlinear operator we have
$$
\|B(u_n,u_n)\|_1 = \sqrt{(D(u_n\cdot \nabla)u_n,D(u_n\cdot \nabla)u_n)}\le c\left( \|Du_n\|_{L^4(\Omega)}^4 + \|u_n\|_{L^{\infty}(\Omega)}^2 |Au_n|^2   \right)^{1/2}
$$
Using the Sobolev inequaities we can write 
$$
\|Du_n\|_{L^4(\Omega)} \le c\|Du_n\|_{L^6(\Omega)}  \le  c|Au_n|
$$
and 
$$
\|u_n\|_{L^{\infty}}  \le  c|Au_n|.
$$
Therefore $\|B(u_n,u_n)\|_1\le c|Au_n|^2$ which implies
$$
f(s)-B(u_n(s),u_n(s)) \in V \quad\mbox{for a.e. $s\in [0,T]$}
$$
and we will have
$$
\|Q_n(B(u_n(s),u_n(s))-f(s))\|_1 \to 0 \quad\mbox{as $n\to \infty$ and for a.e. $s\in [0,T]$}
$$
and
$$
\|Q_n(B(u_n(s),u_n(s))-f(s))\|_1 \le \|B(u_n(s),u_n(s))-f(s)\| \quad\mbox{for a.e. $s\in [0,T]$}.
$$
By the Lebesgue dominated convergence theorem we conclude
$$
\int_0^T\|Q_n(B(u_n(s),u_n(s))-f(s))\|_1\ud s \to 0 \quad\mbox{as $n\to \infty$}.
$$

\end{proof}
A similar result holds for the strong solutions with second order regularity in a periodic domain.%
\begin{theorem}\label{postt2}
$\mathbf{i)}$ Consider the Navier-Stokes equations 
\begin{equation}\label{post2_nse}
\frac{\ud u}{\ud t}+\nu Au+B(u,u)=f\quad\mbox{with}\quad u(0)=u_0
\end{equation}
with $u_0\in V^2$, $f\in L^2(0,T;V^1)\cap L^1(0,T;V^2)$, that hold in a periodic domain. \\
\indent Let $v\in L^{\infty}(0,T;V^2)\cap L^2(0,T;V^3)$ be a numerical approximation of $u$ satisfying
$$
\frac{\ud v}{\ud t}+\nu Av+B(v,v) \in L^1(0,T;V^2)\cap L^2(0,T;V)
$$
and
\begin{eqnarray}
&\|v(0)-u_0\|_2^2 + \int_0^T \|\frac{dv(s)}{ds}(s)+\nu Av(s)
  +B(v(s),v(s))-f(s)\|_2\,\ud s \nonumber\\
&\qquad < \frac{1}{c}\sqrt{\frac{2\nu}{T}}\exp\left(-\int_0^T (c+c')\|u_n\|_3\, \ud t\right).\label{posteriori2}
\end{eqnarray}
Then the Navier-Stokes equations (\ref{post2_nse}) have a strong solution
 $u\in L^{\infty}(0,T;V^2)\cap L^2(0,T;V^3)$.\\ 
\\
\indent $\mathbf{ii)}$ Let $u$ be a strong solution of (\ref{post2_nse}). Then there exists an $N$ such that the Galerkin approximations $u_n$ satisfy the inequality (\ref{posteriori2}) for all $n>N$. Therefore in view of part $\mathbf{(i)}$, $u$ passes \emph{a posteriori} test as a solution approximated by $u_n$ which means that the existence of the strong solution $u$ can be verified using the Galerkin approximations.
\end{theorem}
\\
\indent Using the robustness result for the strong solutions with second order regularity, theorem \ref{robust2}, the proof of the above theorem closely parallels that of theorem \ref{postt}.
%
%%%%%%%%%%%%%%%%%%%%%%%%%%%%%%%%%%%
%        Channel Flow
%%%%%%%%%%%%%%%%%%%%%%%%%%%%%%%%%%%
%
\section{A physical application: channel flow} We considered the Navier-Stokes equations for a flow in a domain with non-moving boundary and with the forcing function in the form of a body force. A physical situation for which the results of this paper and in particular the \emph{a posteriori} test of theorem \ref{postt} can be applied is a channel flow in the domain $0<x<L_x$, $0<y<1$, $0<z<L_z$ with the velocity field $u(x,y,z,t)$ considered to be periodic in $x$ and $z$ with periods $L_x$ and $L_z$ respectively and zero at $y=0$ and $y=1$. The body force $g$ is assumed constant and equal to $(1,1,1)$.
Here we show how to write the inequality (\ref{posteriori}) for the Galerkin approximations, $u_n$, of the solution $u$ of this problem.\\
\indent By theorem \ref{postt} we know that if for some $n$ the Galerkin approximation, $u_n$, satisfies 
\begin{eqnarray*}
&&|Du_n(0)-Du_0| + \int_0^T \|\frac{du_n(s)}{ds}+\nu Au_n(s)+B(u_n(s),u_n(s))
                                                                           -f(s)\|_1\,\ud s \\
&&\quad < \frac{1}{k}\left(\frac{\nu^3}{27T}\right)^{1/4}
                      \exp\left(-\frac{k^2}{2}\int_0^T\frac{27k^2}{2}\frac{1}{\nu^3}|Du_n(s)|^4
                      +\frac{1}{\nu}|Du_n(s)||Au_n(s)|\,\ud s\right)
\end{eqnarray*}
then a strong solution $u$ exists and the Galerkin approximations converge to it. To check this inequality for a Galerkin approximation of the above channel flow example we need to compute  $Au_n$, $B(u_n,u_n)$ and $f=\Pi g$. \\
\indent The functions
$$
w_{\bf{k}}=\e^{2\pi i(\frac{k_1}{L_x}x+\frac{k_3}{L_z}z)}\sin({\pi k_2y})
$$
with ${\bf{k}}=(k_1,k_2,k_3)$, form an orthogonal basis for the space of $L^2$-functions on the cubic domain introduced above, which have periodic bounday values in $x$ and $z$ directions and are zero on $y=0$ and $y=1$. We can therefore define the Galerkin approximation $u_n$ in terms of $w_{\bf k}$ as
\begin{equation}\label{u_expansion}
u_n=\sum_{{\bf k}={\bf n}_0}^{\bf n} \alpha_{\bf{k}}w_{\bf{k}},\quad\mbox{with }\quad\nabla\cdot u_n=0
\end{equation}
where $\alpha_{\bf{k}}=(\alpha_{1\bf{k}}(t),\alpha_{2\bf{k}}(t),\alpha_{3\bf{k}}(t))$, $\alpha_{(-k_1,k_2,-k_3)}=\bar{\alpha}_{(k_1,k_2,k_3)}$, ${\bf n}_0=(-n,0,-n)$ and ${\bf n}=(n,n,n)$. Since $u_n$ is divergence-free we have 
$$
\sum_{{\bf k}={\bf n}_0}^{\bf n} 2\pi i(\frac{k_1}{L_x}\alpha_{1\bf{k}}
                                                     +\frac{k_3}{L_z}\alpha_{3\bf{k}})\sin(\pi k_2y)
                                                     +\pi k_2 \;\alpha_{2\bf{k}}\cos(\pi k_2 y)=0.
$$
After expanding $\cos(\pi k_2 y)$ in terms of sine functions we obtain
$$
\pi i (\frac{k_1}{L_x}\alpha_{1\bf{k}}+\frac{k_3}{L_z}\alpha_{3{\bf k}})+
   \sum_{l=\lfloor\frac{k_2}{2}\rfloor}^{\lfloor\frac{k_2+n}{2}\rfloor}
                (2l+1-k_2)(\frac{1}{2k_2-2l-1}+\frac{1}{2l+1})\alpha_{2\bf{k}}=0.
$$
\indent Defining
$$
\hat{\bf k}=(\frac{\pi i k_1}{L_x},
               \sum_{l=\lfloor\frac{k_2}{2}\rfloor}^{\lfloor\frac{k_2+n}{2}\rfloor}
                (2l+1-k_2)(\frac{1}{2k_2-2l-1}+\frac{1}{2l+1}),\frac{\pi i k_3}{L_z} )
$$
and having (\ref{u_expansion}) we obtain 
\begin{eqnarray*}
Au_n&=& 4\pi^2\sum_{{\bf k}={\bf n}_0}^{\bf n}(\frac{k_1^2}{L_x^2}+\frac{k_2^2}{4}+\frac{k_3^2}{L_z^2})\alpha_{\bf{k}}\;w_{\bf{k}}
\end{eqnarray*}
and 
\begin{eqnarray*}
&&B(u_n,u_n)=2\sum_{{\bf k}=2{\bf n}_0}^{2\bf n}\: \sum_{{\bf j}={\bf k}_{-}} ^{\bf{k}_{+}}
                         \:\sum_{m=\lfloor \frac{k_2+1}{2} \rfloor}^{\lfloor \frac{k_2+2n}{2} \rfloor} 
                        \Bigg\{ (\alpha_{  {\bf k}_m-{\bf j} } \cdot{\bf{j}}) \left( \alpha_{\bf{j}}
                   -\frac{\hat{\bf k}}{|\hat{\bf k}|^2}(\alpha_{\bf{j}}\cdot \hat{\bf{k}}) \right)\\
          &&\qquad    \frac{1}{\pi} \left( \frac{1}{2m-2j_2+1}+\frac{1}{2m-2j_2-2k_2+1}-\frac{1}{2m+1} 
                                          -\frac{1}{2m-2k_2+1} \right) w_{\bf{k}} \Bigg\}
\end{eqnarray*}
where ${\bf k}_m=(k_1,2m+1-k_2,k_3)$, ${\bf k}_{-}=(\min\{k_1,0\},\min\{k_2,0\},\min\{k_3,0\})$ and
${\bf k}_{+}=(\max\{k_1,0\},\max\{k_2,0\},\max\{k_3,0\})$. \\
\indent Now we note that $f=(1,0,1)$, since $f-g=(0,1,0)$ is perpendicular (with respect to $L^2$-norm) to all functions in $H$. Therefore using the expansion of $u_n$ and $Au_n$ with respect to $\{w_{\bf k}\}_{{\bf k}={\bf n}_0}^{\bf n}$ and $B(u_n,u_n)$ with respect to $\{w_{\bf k}\}_{{\bf k}={2\bf n}_0}^{2\bf n}$ we can write (noting that the coefficients of $w_{\bf k}$ when $k_1,k_2$ or $k_3$ are less than $-n$ or bigger than $n$ are zero for $u_n$ and $Au_n$)
$$
E_n=\frac{\ud u_n}{\ud t}+\nu Au_n+B(u_n,u_n)-f=\sum_{{\bf k}=2{\bf n}_0}^{2\bf n} \beta_{\bf k}w_{\bf k}-(1,0,1)
$$
and therefore 
$$
DE_n=\frac{\partial E_n}{\partial x}+\frac{\partial E_n}{\partial y}+\frac{\partial E_n}{\partial z}=2\sum_{{\bf k}=2{\bf n}_0}^{2\bf n} \beta_{\bf k} (\hat{k}_1+\hat{k}_2+\hat{k}_3)w_{\bf k}
$$
and
$$
\|\frac{\ud u_n}{\ud t}+\nu Au_n+B(u_n,u_n)-f\|_1^2=2L_xL_z\sum_{{\bf k}=2{\bf n}_0}^{2\bf n}(\beta_{1{\bf k}}^2+\beta_{2{\bf k}}^2+\beta_{3{\bf k}}^2)(\hat{k}_1+\hat{k}_2+\hat{k}_3)^2.
$$
The norms of $Du_n$ and $Au_n$ are computed as
\begin{eqnarray*}
|Du_n|^2&=& 2L_xL_z\sum_{{\bf k}={\bf n}_0}^{\bf n}(\alpha_{1{\bf k}}^2+\alpha_{2{\bf k}}^2+\alpha_{3{\bf k}}^2)(\hat{k}_1+\hat{k}_2+\hat{k}_3)^2,\\
|Au_n|^2&=& 2\pi^2 L_xL_z \sum_{{\bf k}={\bf n}_0}^{\bf n} (\alpha_{1{\bf k}}^2+\alpha_{2{\bf k}}^2+\alpha_{3{\bf k}}^2)(\frac{k_1^2}{L_x^2}+\frac{k_2^2}{4}+\frac{k_3^2}{L_z^2})^2 .
\end{eqnarray*}
The remaining term on the left-hand side of (\ref{posteriori}), $|Du_n(0)-Du_0|$, is obtained in a similar way to $|Du_n|$.\\
\indent Finally we need the value of the constant $k$ in the inequality (\ref{posteriori}). The cubic domain here also has the strong local Lipschitz property and as it was shown in section 3, for such a domain $k=72(2^{3/4})$.\\

%---------------------
%
\section{Conclusion}

We extended the results of Chernyshenko et al (2006) to the less regular strong solutions not considered in their paper. Unlike their results, the robustness conditions on the data that we obtained for the minimally and second order regular strong solutions in theorems \ref{robust1} and \ref{robust2} depended explicitly on the viscosity coefficient in a way that as the viscosity coefficient became smaller, the conditions (\ref{m1cond}) and (\ref{m2cond}) became more restrictive which obviously is not desirable. 
Moreover, since this is due to the less elegant inequalities available for the trilinear form in the less regular spaces considered here, it does not seem that we can remove the viscosity coefficient from the robustness conditions.\\
\indent On the other hand, the results proved here for the minimally regular strong solutions are valid in a sufficiently smooth bounded domain as well as in the absence of boundaries while in the case of the second order regular solutions considered also here and the more regular ones studied in the above mentioned paper the results are restricted to a periodic domain or the whole of $\mathbb{R}^3$. 
The reason is that the inequalities we need to use in these cases for $|(B(u,v),v )_m|$ 
when $m\ge 2$ are valid in the absence of boundaries where we know that swapping the action of the Leray projector and the Laplacian operator is possible.\\
\indent It would be interesting to see if there exists a big enough $m$ for which it is possible to swap the action of the Laplacian operator and the Leray projector in $D(A^{m/2})$ over a bounded domain. 
If such a finite $m$ exists we can extend the robustness results for regular enough solutions to the more physically relevant bounded domains. 
Furthermore, for two dimensional bounded domains equality of the Stokes and Laplacian operators in $D(A^{m/2})$ can be used in obtaining a better bound on the attractor dimension in $L^2$.\\

\end{document}